\documentclass[a4paper,11pt]{amsart}
 \usepackage[arrow,matrix]{xy}
\usepackage{amsmath,amssymb,amscd,bbm,amsthm,mathrsfs,dsfont,chemarrow}
\usepackage{graphicx}
 \numberwithin{equation}{section}
\newtheorem{thm}{Theorem}[section]
\newtheorem{lem}{Lemma}[section]
\newtheorem{cor}{Corollary}[section]

\theoremstyle{definition}

\theoremstyle{remark} \theoremstyle{example}

\allowdisplaybreaks
\usepackage{geometry}
\begin{document}\numberwithin{equation}{section}
\title[on Matsumoto metrics of scale flag curvature]{on Matsumoto metrics of scale flag curvature \thanks{Supported partially by NNSFC(No. 11171297).}}
\author{Xiaoling Zhang}
\maketitle
\begin{center}
\small {College of Mathematics and Systems Science, Xinjiang University, Urumqi, Xinjiang Province, 830046, P.R.China, xlzhang@ymail.com}
\end{center}
\begin{abstract}
This paper contributes to the study of the Matsumoto metric $F=\frac{\alpha^2}{\alpha-\beta}$, where the $\alpha$ is a Riemannian metric and the $\beta$ is a one form. It is shown that such a Matsumoto metric $F$ is of scalar flag curvature if and only if $F$ is projectively flat.
\end{abstract}
\section{Introduction}

In Finsler geometry, there are several important geometric quantities. In this paper, our main focus is on the flag curvature.

For a Finsler manifold $(M,F)$, the flag curvature $K$ at a point $x$ is a function of tangent planes $P\subseteq T_xM$ and nonzero vectors $y\in P$. This quantity tells us how curved the space is. When $F$ is Riemannian, $K$ depends only on the tangent plane $P\subseteq T_xM$ and is just the sectional curvature in Riemannian geometry. A Finsler metric $F$ is said to be of scalar flag curvature if the flag curvature $K$ at a point $x$ is independent of the tangent plane $P\subseteq T_xM$, that is, the flag curvature $K$ is a scalar function on the slit tangent bundle $TM\backslash \{0\}$. It is known that every locally projectively flat Finsler metric is of scalar flag curvature. However, the converse is not true.

$(\alpha, \beta)$-metrics form a special and important class of Finsler metrics which can be expressed in the form $F=\alpha\phi(s)$, $s=\frac {\beta}{\alpha}$, where $\alpha=\sqrt{a_{ij}(x)y^iy^j}$ is a Riemannian metric, $\beta=b_i(x)y^i$ is a 1-form on $M$ and $\phi(s)$ is a $C^{\infty}$ positive function satisfying \eqref{ta8} on some open interval. In particular, when $\phi(s)=\frac 1{1-s}$, the Finsler metric $F=\frac{\alpha^2}{\alpha-\beta}$ is called a Matsumoto metric, which was first introduced by M.Matsumoto to study the time it takes to
negotiate any given path on a hillside(cf. \cite{AIM}). Recently, some researchers have studied Matsumoto metrics(\cite{RR1,RR2,XZ,ZS2}).

Randers metrics are the simplest $(\alpha, \beta)$-metrics. Bao, etc., finally classify Randers metrics of constant flag curvature by using the navigation method (see \cite{BRS}). Further, Shen, etc., classify Randers metrics of weakly isotropic flag curvature (see \cite{SY}). There are Randers metrics of scalar flag curvature which are not of weakly isotropic flag curvature or not locally projectively flat(see \cite{BZ,SXql}). Besides, some relevant researches are refereed to \cite{CxyS,SXh,Y1}, under additional conditions. So far, Randers metrics of scalar flag curvature are still unknown. Xia and Yang obtain the results for Kropina metrics and m-Kropina metrics, respectively(see \cite{X,Y2}).

Our main result concerns Matsumoto metrics of scalar flag curvature.

\begin{thm}\label{tta1}
Let $F=\frac{\alpha^2}{\alpha-\beta}$ be a non-Riemannian Matsumoto metric on an $n$-dimensional manifold $M$, $n\geq3$. Then $F$ is of scalar flag curvature if and only if $F$ is projectively flat, i.e., $\alpha$ is locally projectively flat and $\beta$ is parallel with respect to $\alpha$.
\end{thm}

Li obtains that an $n(\geq3)$-dimensional non-Riemannian Matsumoto metric is projectively flat if and only if $\alpha$ is locally projectively flat and $\beta$ is parallel with respect to $\alpha$(see \cite{L}). It is known that $\alpha$ is locally projectively flat is equivalent to that $\alpha$ is of constant curvature. Hence, a Matsumoto metric, which is projectively flat (i.e., of scalar flag curvature), must be locally Minkowskian.

The content of this paper is arranged as follows. In \S \ref{t2} we introduce essential curvatures of Finsler metrics, as well as notations and conventions. And we give basic formulas for proofs of Theorems in the following section. In \S \ref{t4} the characterization of scalar flag curvature is given under the assumption that the dual of $\beta$, with respect to $\alpha$, is a constant Killing vector field. By using it, Theorem \ref{tta1} is proved in \S \ref{t5}.

\section{Preliminaries}\label{t2}
In this section, we give a brief description of several geometric quantities in Finsler geometry. See \cite {CS} in detail.

Let $F$ be a Finsler metric on an $n$-dimensional smooth manifold $M$ and $(x, y)=(x^i, y^i)$ the local coordinates on the tangent bundle $TM$. The geodesics of F are locally characterized by a system of second order ordinary differential equations
\begin{eqnarray*}
\frac{d^2x^i}{dt^2}+2G^i\left(x(t), \frac{dx(t)}{dt}\right)=0,  \label{ta1}\end{eqnarray*}
where
\begin{eqnarray*}G^i=\frac 14 g^{ij}\{[F^2]_{x^ky^j}y^k-[F^2]_{x^j}\}.\label{ta2}\end{eqnarray*}
$G^i$ are called spray coefficients of F.

The Riemann curvature $R$ (in local coordinates ${R^i}_k\frac{\partial}{\partial x^i}\otimes dx^k$) of $F$ is defined by
\begin{eqnarray*}
{R^i}_k=2\frac{\partial G^i}{\partial x^k}-\frac{\partial^2G^i}{\partial x^j\partial y^k}y^j
+2G^j\frac{\partial^2 G^i}{\partial y^j\partial y^k}-\frac{\partial G^i}{\partial y^j}\frac{\partial G^j}{\partial y^k}.\label{ta3}
\end{eqnarray*}
It is known that $F$ is of scalar flag curvature if and only if, in a local coordinate system,
\begin{eqnarray*}
{R^i}_k=K(x,y)\{F^2\delta^i_{\,\,\,k}-FF_{y^k}y^i\}.\label{ta4}
\end{eqnarray*}
In particular,
$F$ is of weakly isotropic(resp.isotropic or constant) flag curvature if $K(x,y)=\frac {3\theta}F+\kappa(x)$(resp. $K(x,y)=\kappa(x)$ or constant), where $\kappa=\kappa(x)$ is a scalar function and $\theta$ is a 1-form on $M$. In dimension $n\geq 3$, a Finsler metric $F$ is of isotropic flag curvature if and only if $F$ is of constant flag curvature by Schur's Lemma. In general, a Finsler metric of weakly isotropic flag curvature and that of isotropic flag curvature are not equivalent.

The Ricci curvature(or Ricci scalar) of $F$ is defined by \begin{eqnarray*}Ric:=R^m_{\,\,\,m}.\label{2.5}\end{eqnarray*}

For a vector $y\in T_xM \setminus \{0\}$, define $W_y(u)=W^i_{\,\,\,j}(y)u^j\frac{\partial}{\partial x^i}|_x$ by
 \begin{eqnarray}
W^i_{\,\,\,j}(y):=A^i_{\,\,\,j}-\frac{1}{n+1}\frac{\partial A^k_{\,\,\,j}}{\partial y^k}y^i,\label{ta5}
\end{eqnarray}
where
\begin{eqnarray}
A^i_{\,\,\,j}:=R^i_{\,\,\,j}-\frac{1}{n-1}R^m_{\,\,\,m}\delta^i_{\,\,\,j}. \label{ta6}
\end{eqnarray}
We call $W:=\{W_y\}_{y\in T_xM \setminus \{0\}}$ the Weyl curvature. It is easy to check that $W$ is a projective invariant which means that if $G^i=\widetilde{G}^i+Py^i$, then $W^i_{\,\,\,j}=\widetilde{W}^i_{\,\,\,j}$. As is known that a Finsler metric is of scalar flag curvature if and only if its Weyl curvature $W=0$. For a Riemannian space $(M, g)$ of dimension $n\geq3$, the following conditions are equivalent: (a) $W=0$, (b) $g$ is of scalar flag curvature, (c) $g$ is of constant curvature, (d) $g$ is locally projectively flat.

($\alpha$, $\beta$)-metrics are an important class of Finsler metrics. An ($\alpha$, $\beta$)-metric is expressed by the following form,
\begin{eqnarray*}
F=\alpha\phi(s), \qquad s:=\frac{\beta}{\alpha},\label{ta7}
\end{eqnarray*}
where $\alpha=\sqrt{a_{ij}(x)y^iy^j}$ is a Riemannian metric and $\beta=b_i(x)y^i$ is a 1-form. $\phi=\phi(s)$ is a $C^{\infty}$ positive function on an open interval $(-b_0, b_0)$ satisfying
\begin{eqnarray}\phi(s)-s\phi'(s)+(b^2-s^2)\phi''(s)>0, \ \ \ |s|\leq b<b_0, \label{ta8}\end{eqnarray}
where $b:=\|\beta(x)\|_\alpha=\sqrt{a^{ij}b_ib_j}$. It is known that $F=\alpha\phi(s)$ is a Finsler metric if and only if $\parallel \beta(x)\parallel_{\alpha}<b_0$ for any $x\in M$. In particular, if $\phi(s)=1+s$, then it is a Randers metric. If $\phi(s)=\frac 1{1-s}$, then it is a Matsumoto metric $F=\frac{\alpha^2}{\alpha-\beta}$, etc. It is known that a Matsumoto metric $F$ is a Finsler metric if and only if $b<b_0=\frac 12$.

Let $G^i(x,y)$ and ${^\alpha G}^i(x,y)$ denote spray coefficients of $F$ and $\alpha$, respectively. To express formulae for spray coefficients $G^i$ of $F$ in terms of $\alpha$ and $\beta$, let's introduce some notations. Let $b_{i|j}$ be a horizontal covariant derivative of $b_i$ with respect to $\alpha$. Denote
\begin{equation*}\begin{aligned}
r_{ij}:=&\frac 12(b_{i|j}+b_{j|i}), \qquad \quad s_{ij}:=\frac 12(b_{i|j}-b_{j|i}),\\
{s^i }_j:=&a^{ik}s_{kj},\quad s_j:=b_i{s^i }_{  j}=b^is_{ij},\quad r_j:=b^ir_{ij},\quad r:=b^ir_i, \\
r_0:=&r_iy^i,\quad s_0:=s_iy^i, \quad r_{00}:=r_{ij}y^iy^j.
\end{aligned}\end{equation*}
Throughout this paper, we use $a_{ij}$ to raise and lower the indices of $b_i$, $r_{ij}$, $s_{ij}$, $r_i$, $s_i$ and $y^i$, etc.  With these, we can express the spray coefficients $G^i$ as follows(see \cite{CST})
 \begin{eqnarray}G^i={^\alpha G}^i+\alpha {Qs^i }_{  0}+\Theta(r_{00}-2\alpha Qs_0)\frac{y^i}{\alpha}+
 \psi(r_{00}-2\alpha Qs_0)b^i,\label{ta9}\end{eqnarray}
where
\begin{eqnarray*}Q:&=&\frac {\phi'}{\phi-s\phi'},\\
\Theta:&=&\frac{\phi'(\phi-s\phi')}{2\phi[\phi-s\phi''+(B-s^2)\phi'']}-s\psi,\\
\psi:&=&\frac{\phi''}{2[\phi-s\phi'+(B-s^2)\phi'']}.\\ \end{eqnarray*}
Here $B:=b^2$. In particular, for a Matsumoto metric $F=\frac {\alpha^2}{\alpha-\beta}$,
it follows from (\ref{ta9}) that
\begin{eqnarray}G^i={^\alpha G}^i+\alpha {Qs^i }_{  0}+\Theta(r_{00}-2\alpha Qs_0)\frac{y^i}{\alpha}+ \psi(r_{00}-2\alpha Qs_0)b^i,\label{ta10}\end{eqnarray}
where
\begin{eqnarray}Q=\frac{1}{A_1}, \,\,\Theta=\frac{1-4s}{2A_2}, \,\, \psi=\frac{1}{A_2},  \label{ta11}\end{eqnarray}
and $A_1=1-2s$, $A_2=1+2B-3s$ are functions of $s, B$ respectively. Note that $B<\frac 14$ for a Matsumoto metric $F=\frac {\alpha^2}{\alpha-\beta}$. The variations $s$ or $B$ is also a linear combination of two functions $A_1$ and $A_2$, which means $A_1, A_2$ can be regarded as new variations in place of the variations $s, B$.

\section{Basic Formulas for $(\alpha,\beta)$-metrics} \label{t3}
From (\ref{ta10}), we have that
\begin{eqnarray}G^i={^\alpha G}^i+Py^i+Q^i,   \label{te1}
\end{eqnarray}
where $P=\frac{1}{\alpha}(r_{00}-2\alpha Qs_0)\Theta, Q^i=\alpha {Qs^i }_{  0}+ \psi(r_{00}-2\alpha Qs_0)b^i$. Let \begin{eqnarray}\overline{G}^i={^\alpha G}^i+Q^i.   \label{te2}\end{eqnarray}
By (\ref{te1}), (\ref{te2}) and (\ref{ta5}), we have
\begin{eqnarray} W^i_{\,\,\,j}=\overline{W}^i_{\,\,\,j}. \label{te5}\end{eqnarray}

Assume that $F$ is of scalar flag curvature. Now we have
\begin{equation}\begin{aligned}  \label{te6}
0&=W^i_{\,\,\,j}=\overline{W}^i_{\,\,\,j}  \\
&=\overline{R}^i_{\,\,\,j}-\frac{\overline{R}^m_{\,\,\,m}}{n-1}\delta^i_{\,\,\,j}
-\frac{1}{n+1}(\overline{R}^m_{\,\,\,j}-\frac{\overline{R}^s_{\,\,\,s}}{n-1}\delta^m_{\,\,\,j})_{.\,m}y^i,
\end{aligned}\end{equation}
Contracting (\ref{te6}) with $b^jb_i$ yields
\begin{eqnarray}
0&=&\overline{R}^i_{\,\,\,j}b^jb_i-B\frac{\overline{R}^m_{\,\,\,m}}{n-1}
-\frac{s\alpha}{n+1}(\overline{R}^i_{\,\,\,j}b^j-\frac{\overline{R}^m_{\,\,\,m}}{n-1}b^i)_{.\,i}.   \label{te7}\end{eqnarray}

Using formulas and notations in \cite{CST}, we have
\begin{equation}\label{tc1}\begin{aligned}
\overline{R}^i_{\,\,\,j}b^j&={^\alpha R}^i_{\,\,\,j}b^j+  b^i(sC_{21}+BC_{22}+rC_{24}-C_{230}s_0+C_{240}r_0-v\alpha b^js_{j|0}+2v\alpha b^js_{0|j}+vQ\alpha^2s_ks^k\\
&+2\psi b^jr_{00|j}-2\psi b^jr_{j0|0}+2Q\psi\alpha s^kr_{0k}+4Q\psi\alpha r^ks_{k0})\\
&+s^i(sC_{31}+BC_{32}+2Q\psi\alpha r_0)+s^i_{\,\,\,0}(sC_{310}+BC_{320}-C_{331}s_0)\\
&+r^iC_4+r^i_{\,\,\,0}(sC_{410}+BC_{420}-2\psi r_0)+s^i_{\,\,\,k}s^k_{\,\,\,0}(sC_{332}+BC_{333})+s^i_{\,\,\,k}s^kQ^2\alpha^2   \\
&+2b^js^i_{\,\,\,0|j}Q\alpha-b^js^i_{\,\,\,j|0}Q\alpha+s^i_{\,\,\,0|0}(sC_{311}-BQ_s),
\end{aligned}\end{equation}
where
\begin{eqnarray}
v(s,B)&:=&-2\psi Q,\nonumber \\
C_{21}&=&\frac{r_{00}^2}{\alpha^2}
 ((s\psi_{s}^2-2s\psi\psi_{ss}-2\psi\psi_{s})(B-s^2)+s\psi_{ss}+\psi_{s}+4s^2\psi\psi_{s})  \nonumber \\
 &&+\frac{r_{00}s_{0}}{\alpha}((2s\psi_{s}v_s-3v\psi_{s}-2sv\psi_{ss}-2s\psi v_{ss})(B-s^2)+sv_{ss}+2s\psi_{Bs}  \nonumber \\
 &&-2sQ\psi_{ss}+5s^2v\psi_{s}+sQ_s\psi_{s}+2s^2v_s\psi-2sv\psi-3Q\psi_{s})  \nonumber \\
 &&+\frac{2r_{00}r_{0}}{\alpha}(-s\psi\psi_{s}+s\psi_{Bs})+s_{0}^2((sv_s^2-2svv_{ss}-vv_s)(B-s^2)+3s^2vv_s  \nonumber \\
 &&-2sQv_{ss}-3sv^2-Qv_s-2v_B+sQ_sv_s+2sv_{Bs})  \nonumber \\
 &&+2s_{0}r_{0}(sv_{Bs}-\psi  v+s\psi v_s-2sv\psi_{s}-v_B)+\frac{r_{00|0}}{\alpha}s\psi_{s}  \nonumber \\
 &&+2r_{k0}s^k_{\,\,\,0}(-Q\psi+s\psi Q_s-2sQ\psi_{s})+s_{0|0}(sv_s-v)  \nonumber \\
 &&+\alpha s_ks^k_{\,\,\,0}(sQ_sv+Qv-2sQv_s),  \nonumber \\
C_{22}&=&\frac{r_{00}^2}{\alpha^2}((2\psi\psi_{ss}-\psi_{s}^2)(B-s^2)-\psi_{ss}-4s\psi\psi_{s})  \nonumber \\
&&+\frac{r_{00}s_0}{\alpha}(2(v\psi_{ss}+\psi v_{ss}-\psi_{s}v_{s})(B-s^2)-2\psi_{Bs}-v_{ss}+2v\psi \nonumber \\
&&-2s\psi v_{s}-5sv\psi_{s}+2Q\psi_{ss}-Q_{s}\psi_{s})  \nonumber \\
&&+\frac{2r_{00}r_0}{\alpha}(-\psi_{Bs}+\psi\psi_{s})+2s_0r_0(-\psi v_{s}+2\psi_ {s}v-vB_s)  \nonumber \\
&&+s^2_0((2vv_{ss}-v^2_{s})(B-s^2)-Q_sv_{s}-2v_{Bs}+2v^2+2Qv_{ss}-3svv_s)  \nonumber \\
&&-\frac{r_{00|0}}{\alpha}\psi_{s}-s_{0|0}v_{s}+as_ks^k_{\,\,\,0}(2Qv_{s}-Q_{s}v)+2r_{k0}s^k_{\,\,\,0}(2Q\psi_{s}-Q_{s}\psi),   \nonumber \\
C_{23}&=&r_{00}((-v\psi_{s}+2\psi v_{s})(B-s^2)+4\psi_B-Q\psi_{s}-v_{s}+2sv\psi)\nonumber \\
 &&+\alpha s_{0}(vv_{s}(B-s^2)+2v_B+Qv_{s}+sv^2)-2\alpha r_{0}(\psi v+v_B),\nonumber \\
C_{24}&=&4r_{00}(\psi^2+\psi_{B})+4\alpha s_0(v_B+\psi v),    \nonumber \\
C_{230}&=&\frac{3r_{00}}{\alpha}(1+s Q)\psi_{s}+3s_0(v_{s}-(v-sv_{s})Q),  \nonumber \\
C_{240}&=&\frac{r_{00}}{\alpha}(2\psi\psi_{s}(B-s^2)-\psi_{s})-4r_0(\psi^2+\psi_{B})  \nonumber \\
&&+s_0((4v\psi_{s}-2\psi v_{s})(B-s^2)+v_{s}-2sv\psi+4Q\psi_{s}-4\psi_{B}),  \nonumber \\
C_{31}&=&r_{00}(-2Q\psi+2sQ_{s}\psi-sQ\psi_{s})-\alpha s_0(vQ+sQv_{s}-2sQ_{s}v),  \nonumber \\
C_{32}&=& r_{00}(-2Q_{s}\psi+Q\psi_{s})-\alpha s_0(2Q_{s}v-Qv_{s}),  \nonumber \\
C_{310}&=&\frac{r_{00}}{\alpha }((sQ_{s}\psi_{s}-2s\psi Q_{ss})(B-s^2) \nonumber \\
&&+2s^2Q_{s}\psi-2sQ\psi+s^2Q\psi_{s}+sQ_{ss}+s\psi_{s})  \nonumber \\
&&+s_0((sQ_{s}v_{s}-Q_{s}v-2sQ_{ss}v)(B-s^2)+s^2Qv_{s}-QQ_{s}+sQ_{s}^2-3sQv  \nonumber \\
&&-v+sv_{s}-2sQQ_{ss}+2s^2Q_{s}v),   \nonumber \\
C_{320}&=&\frac{r_{00}}{\alpha}((2Q_{ss}\psi-Q_s\psi_{s})(B-s^2)+2Q\psi-2s\psi  Q_s-sQ\psi_{s}-Q_{ss}-\psi_{s})\nonumber \\ &&+s_{0}((2Q_{ss}v-Q_sv_{s})(B-s^2)-v_{s}-Q^2_s+2Qv-2sQ_sv+2QQ_{ss}-sQv_{s}),  \nonumber \\
C_{331}&=&-3Q^2+3sQQ_s+3Q_s,  \nonumber \\
C_{4}&=&2\psi r_{00}+2v\alpha s_{0},  \nonumber \\
C_{410}&=&\frac{r_{00}}{\alpha}s\psi_{s}+s_{0}(sv_{s}-v),  \nonumber \\
C_{420}&=&-\frac{r_{00}}{\alpha}\psi_{s}-s_{0}v_{s},  \nonumber \\
C_{332}&=&\alpha(Q-sQ_s)Q,  \nonumber \\
C_{333}&=&\alpha QQ_s,  \nonumber \\
C_{311}&=&sQ_s-Q,  \label{tc2}\end{eqnarray}
and
\begin{eqnarray}\overline{R}^m_{\,\,\,m}&=&{^\alpha R}^m_{\,\,\,m}
+\frac{r_{00}^2}{\alpha^2}c_{2}+\frac{1}{\alpha}(r_{00}s_0c_4+r_{00}r_0c_6+r_{00|0}c_8)\nonumber\\
&&+s^2_0c_{10}+(rr_{00}-r^2_{0})c_{11}+r_0s_0c_{13}+(r_{00}r^m_{\,\,\,m}-r_{0m}r^m_{\,\,\,0}+r_{00|m}b^m-r_{0m|0}b^m)c_{14} \nonumber\\
&&+r_{0m}s^m_{\,\,\,0}c_{16}+s_{0|0}c_{18}+s_{0m}s^m_{\,\,\,0}c_{19}\nonumber\\
&&+\alpha\{rs_0c_{20}+s_{m}s^m_{\,\,\,0}c_{22}+(3s_{m}r^m_{\,\,\,0}-2s_0r^m_{\,\,\,m}+2r_ms^m_{\,\,\,0}-2s_{0|m}b^m+s_{m|0}b^m)c_{23}+s^m_{\,\,\,0|m}c_{24}\}  \nonumber\\
&&+\alpha^2(s_{m}s^mc_{25}+s^i_{\,\,\,m}s^m_{\,\,\,i}c_{26}),   \label{td1}\end{eqnarray}
where
\begin{eqnarray}
c_2&=&(2 \psi \psi_{ss}- \psi_{s}^2) ( B-s^2)^2- (6 s \psi \psi_{s} +\psi_{ss}) ( B-s^2)+2 s \psi_{s},  \nonumber\\
c_4&=&(- 4 \psi (2 Q \psi_{ss}+ Q_{s} \psi_{s} + Q_{ss} \psi) + 4Q \psi_{s}^2) ( B - s^2)^2  \nonumber\\
&&+(- 4 \psi^2 ( Q - s Q_{s})+2 (2 Q_{ss} \psi+Q_{s} \psi_{s}+2 Q \psi_{ss})-2 \psi_{sB}+20 s Q \psi \psi_{s}) ( B-s^2)  \nonumber\\
&&+ 2 \psi (Q - s Q_{s}) - 4 \psi_{s} - Q_{ss}- 10 s Q \psi_{s}, \nonumber\\
c_6&=&2 (2 \psi \psi_{s} - \psi_{sB}) ( B-s^2) - 2\psi_{s},  \nonumber\\
c_8&=&-\psi_{s} (B-s^2), \nonumber\\
c_{10}&=&(4 \psi^2 (2 Q Q_{ss}- Q_{s}^2)+ 8 Q \psi (Q \psi_{ss}+ Q_{s} \psi_{s}) - 4 Q^2 \psi_{s}^2) (B-s^2)^2  \nonumber\\
&&+(- 16 s Q \psi (Q \psi_{s} + Q_{s} \psi)- 4 \psi (2 Q Q_{ss}- Q_{s}^2) - 4 Q ( Q \psi_{ss}+ Q_{s} \psi_{s})  \nonumber\\
&&+4 (Q \psi_{sB} + Q_{s} \psi_B)+8Q^2\psi^2) (B-s^2)  \nonumber\\
&&- 4 s^2 Q^2 \psi^2+4 (2+ 3 s Q) ( Q \psi_{s} + Q_{s} \psi)-8 Q^2 \psi+2 Q Q_{ss}- Q_{s}^2+4 s Q \psi_{B},  \nonumber\\
c_{11}&=&4 \psi^2+4 \psi_{B}, \nonumber\\
c_{13}&=&(8 \psi ( Q_{s} \psi - Q \psi_{s})+ 4 (Q \psi_{sB} + Q_{s} \psi_B)) (B-s^2)+ 8 s Q \psi^2+4 Q \psi_{s}-4(1-s Q) \psi_{B}, \nonumber\\
c_{14}&=&2 \psi, \nonumber\\
c_{16}&=&-4 ( Q_{s} \psi - Q \psi_{s}) (B-s^2)+ 2 Q_{s} - 2 (1 + 2 s Q) \psi,  \nonumber\\
c_{18}&=&2 (Q_{s} \psi+Q \psi_{s}) (B-s^2) - Q_{s}+2 s Q \psi,\nonumber\\
c_{19}&=&-2 Q^2+2(1+s Q) Q_{s},\nonumber\\
c_{20}&=&-8 Q (\psi^2+\psi_{B}),  \nonumber\\
c_{22}&=&-4 Q^2 \psi_{s} (B-s^2)+ 2 Q \psi, \nonumber\\
c_{23}&=&2 Q \psi,  \nonumber\\
c_{24}&=&2 Q,  \nonumber\\
c_{25}&=&-4 Q^2 \psi,  \nonumber\\
c_{26}&=&-Q^2.
\label{td2}\end{eqnarray}

\section{A Special Case for Matsumoto metrics}\label{t4}
\begin{lem}\label{s} (\cite{XZ}) (1) $1-s^2$ and $(1+2B)^2-9s^2$ are relatively prime polynomials in $s$ if and only if $B\neq 1$;\,
(2) $B-s^2$ and $(1+2B)^2-9s^2$ are relatively prime polynomials in $s$ if and only if $B\neq 1$;\,(3) $B-s^2$(or $(1+2B)^2-9s^2$) and $1-4s^2$ are relatively prime polynomials in $s$ if and only if $B\neq \frac 14$.
\end{lem}

For an $(\alpha,\beta)$-metrics, the form $\beta$ is said to be Killing (resp. closed) form if $r_{ij}=0$\,\,(resp. $s_{ij}=0$). $\beta$ is said to be a constant Killing form if its dual is a Killing vector field and has constant length with respect to $\alpha$, equivalently $r_{ij}=0,s_i=0$.

\begin{lem}\label{ttq1}
Let $F=\frac{\alpha^2}{\alpha-\beta}$ be a non-Riemannian Matsumoto metric on an $n$-dimensional manifold $M$, $n\geq3$. Suppose $\beta$ is a constant Killing form. Then $F$ is of scalar flag curvature if and only if $F$ is projectively flat.
\end{lem}
\noindent{\it Proof.} Suppose $\beta$ is a constant Killing form, i.e., $r_{ij}=0,s_i=0$. We need to rewrite the equation \eqref{te7}. By \eqref{tc1}, we have
\begin{equation}\label{tq1}\begin{aligned}
\overline{R}^i_{\,\,\,j}b^j&={^\alpha R}^i_{\,\,\,j}b^j+\alpha [(B-s^2)QQ_s+sQ^2+Q]s^i_{\,\,\,k}s^k_{\,\,\,0}+2\alpha Q b^js^i_{\,\,\,0|j}-[(B-s^2)Q_s+sQ]s^i_{\,\,\,0|0},
\end{aligned}\end{equation}

(1) The calculations for $\overline{R}^i_{\,\,\,k}b^kb_i$.

By \eqref{tq1}, we have
\begin{equation}\label{tq2}\begin{aligned}
\overline{R}^i_{\,\,\,j}b^jb_i=&{^\alpha R}^i_{\,\,\,j}b^jb_i-[(B-s^2)Q_s+sQ]s_{0k}s^k_{\,\,\,0}.
\end{aligned}\end{equation}

(2) The calculations for ${\overline{R}}^m_{\,\,\,m}$.

By \eqref{td1}, ${\overline{R}}^m_{\,\,\,m}$ can be rewritten as
\begin{equation}\label{tq3}\begin{aligned}
{\overline{R}}^m_{\,\,\,m}=&{^\alpha R}^m_{\,\,\,m}+2(-Q^2+sQQ_s+Q_s)s_{0k}s^k_{\,\,\,0}
+2\alpha Qs^k_{\,\,\,0|k}-\alpha^2Q^2s^i_{\,\,\,k}s^k_{\,\,\,i}.
\end{aligned}\end{equation}

(3) The calculations for $(\overline{R}^i_{\,\,\,j}b^j)_{.\,i}$.

By \eqref{tq1}, $(\overline{R}^i_{\,\,\,j}b^j)_{.\,i}$ can be rewritten as
\begin{equation}\label{tq4}\begin{aligned}
(\overline{R}^i_{\,\,\,j}b^j)_{.\,i}=&({^\alpha R}^i_{\,\,\,j})_{.\,i}
+\frac{1}{\alpha}(B-s^2)(QQ_s-sQ^2_s-sQQ_{ss}-Q_{ss})s_{0k}s^k_{\,\,\,0}\\
&-[(B-s^2)Q_s+sQ]s^k_{\,\,\,0|k}  +\alpha[(B-s^2)QQ_s+sQ^2+Q]s^i_{\,\,\,k}s^k_{\,\,\,i}.
\end{aligned}\end{equation}

(4) The calculations for $(\overline{R}^m_{\,\,\,m}b^i)_{.\,i}$.

By \eqref{tq3} and a direct computation, we obtain
\begin{equation}\label{tq5}\begin{aligned}
(\overline{R}^m_{\,\,\,m})_{.\,i}b^i=&({^\alpha  R}^m_{\,\,\,m})_{.\,i}b^i +\frac{2}{\alpha}(B-s^2)(-QQ_s+sQ^2_s+sQQ_{ss}+Q_{ss})s_{0k}s^k_{\,\,\,0}\\
&+2[(B-s^2)Q_s+sQ]s^k_{\,\,\,0|k}  -2\alpha[(B-s^2)QQ_s+sQ^2+Q]s^i_{\,\,\,k}s^k_{\,\,\,i}.
\end{aligned}\end{equation}

Notice that $A_1^{-4}=(1-2s)^{-4}$ only exist in the term $(\overline{R}^i_{\,\,\,j}b^j)_{.\,i}$ and $(\overline{R}^m_{\,\,\,m}b^i)_{.\,i}$. Let \begin{equation*}\label{tq6}\begin{aligned}
V:=&\frac{1}{\alpha}(B-s^2)(sQ^2_s+sQQ_{ss})s_{0k}s^k_{\,\,\,0}\\
&=\frac{12s(B-s^2)}{\alpha A_1^4}s_{0k}s^k_{\,\,\,0}\end{aligned}\end{equation*}
Plugging \eqref{tq2}-\eqref{tq5} into \eqref{te7} yields
\begin{equation}\label{tq7}\begin{aligned}
0=&-\frac{s\alpha}{n+1}(-V-\frac{2}{n-1}V)+(\cdots)\\
=&\frac{12s^2(B-s^2)}{(n-1)A_1^4}s_{0k}s^k_{\,\,\,0}+(\cdots),\\
=&\frac{12s^2(B-s^2)(1+4s)^4}{(n-1)(1-4s^2)^4}s_{0k}s^k_{\,\,\,0}+(\cdots),\\
 \end{aligned}\end{equation}
where $(\cdots)$ does not contain the factor $A_1^{-4}$.

For the Matsumoto metric $F$, we have $B<1/4$, which implies $B-s^2$ can not be divided by $1-4s^2$ from Lemma \ref{s}. Obviously, $s^2$ can not be divided by $1-4s^2$ either. Thus $s_{0k}s^k_{\,\,\,0}$ must be divided by $1-4s^2=\frac{\alpha^2-4\beta^2}{\alpha^2}$. This implies there exists some function $\rho=\rho(x)$, such that $s_{0k}s^k_{\,\,\,0}=\rho(\alpha^2-4\beta^2)$. Differentiating both sides of it with $y^iy^j$ yields $s_{kj}s^k_{\,\,\,i}=\rho(a_{ij}-4b_ib_j)$. Then contracting it with $b^j$ gives $0=\rho(1-4B)b_j$, which means $\rho=0$ and $0=s_{kj}s^k_{\,\,\,i}$. Since $|s_{ij}|^2_{\alpha}=-s^i_{\,\,\,k}s^k_{\,\,\,i}=-a^{ij}s_{jk}s^k_{\,\,\,i}=0$, $s_{ij}=0$. Hence $b_{i|j}=0$, which means that $\beta$ is parallel with respect to $\alpha$. So by \eqref{tq7}, we have $W^i_{\,\,\,j}=\overline{W}^i_{\,\,\,j}={^\alpha W}^i_{\,\,\,j}=0$. Then $\alpha$ is projectively flat, so is $F$.

The converse is obvious. This completes the proof of Lemma \ref{ttq1}. \hfill$\Box$

\section{Proof of Theorem \ref{tta1}} \label{t5}
\begin{lem}\label{conf} Let $F=\frac {\alpha^2}{\alpha-\beta}$ be a non-Riemannian Matsumoto metric on an $n(\geq2)$-dimensional manifold $M$. If $F$ is of scalar flag curvature, then $\beta$ is a conformal 1-form with respect to $\alpha$. i.e. there is a function $\sigma=\sigma(x)$ on $M$ such that $r_{00}=\sigma\alpha^2$.
\end{lem}
\noindent{\it Proof.} Rewriting the equation \eqref{te7}.

(1) The calculations for $\overline{R}^i_{\,\,\,j}b^jb_i$.

By \eqref{tc1}, we have
\begin{equation}\label{tp1}\begin{aligned}
\overline{R}^i_{\,\,\,j}b^jb_i&={^\alpha R}^i_{\,\,\,j}b^jb_i+  B(sC_{21}+BC_{22}+rC_{24}-C_{230}s_0+C_{240}r_0-\alpha v b^js_{j|0}+2\alpha v b^js_{0|j}+\alpha^2vQs_ks^k\\
&+2\psi b^jr_{00|j}-2\psi b^jr_{j0|0}+2\alpha Q\psi s^kr_{0k}+4\alpha Q\psi r^ks_{k0})\\
&  +s_{0}(sC_{310}+BC_{320}-C_{331}s_0)+rC_4+r_{0}(sC_{410}+BC_{420}-2\psi r_0)+s_{k}s^k_{\,\,\,0}(sC_{332}+BC_{333})\\
&+\alpha^2Q^2 s_ks^k +2\alpha Qb_ib^js^i_{\,\,\,0|j}-\alpha Qb_ib^js^i_{\,\,\,j|0}+(sC_{311}-BQ_s)b_is^i_{\,\,\,0|0},
\end{aligned}\end{equation}
which can be rewritten as
\begin{equation}\label{tp2}\begin{aligned}
\overline{R}^i_{\,\,\,j}b^jb_i&={^\alpha R}^i_{\,\,\,j}b^jb_i+\frac{f_1(s,B)}{\alpha^{2}A^4_2}r^2_{00}+\sum_{i_1,j_1,k_1}\frac 1{\alpha^{i_1}A_1^{j_1}A_2^{k_1}}(\cdots),
\end{aligned}\end{equation}
where the sum indices $i_1,j_1$ and $k_1$ satisfy $-2\leq i_1\leq 1, 0\leq j_1\leq 3$ and $0\leq k_1\leq 4$ respectively, $f_1(s,B)$ is a nonzero polynomial of $s, B$ given by
\begin{equation*}\label{tp3}\begin{aligned}
&f_1(s,B)/A^4_2=B(B-s^2)^2(2\psi\psi_{ss}-\psi_{s}^2)- B(B-s^2)(6s\psi\psi_{s}+\psi_{ss})+Bs\psi_{s},\\
&f_1(s,B)=-3B\{9s^4+3(1-4B)s^3+(8B^2-16B-1)s+3B(B+2)\},
\end{aligned}\end{equation*}
$(\cdots)$, independent of $r_{00}^2$, is a polynomial of $A_1, A_2$(or $s, B$), and the degree of $(\cdots)$ in $s$ is no more than deg$(A_1^{j_1}A_2^{k_1})=j_1+k_1$.\\

(2) The calculations for ${\overline{R}}^m_{\,\,\,m}$.

By \eqref{td1}, $B{\overline{R}}^m_{\,\,\,m}$ can be rewritten as
\begin{equation}\label{tp4}\begin{aligned}
B{\overline{R}}^m_{\,\,\,m}&=B{^\alpha R}^m_{\,\,\,m}+\frac{f_2(s,B)}{\alpha^{2}A^4_2}r^2_{00}+\sum_{i_2,j_2,k_2}\frac 1{\alpha^{i_2}A_1^{j_2}A_2^{k_2}}(\cdots),
\end{aligned}\end{equation}
where the sum indices $i_2,j_2$ and $k_2$ satisfy $-2\leq i_2\leq 1, 0\leq j_2\leq 3$ and $0\leq k_2\leq 4$ respectively, $f_2(s,B)$ is a nonzero polynomial of $s, B$ given by
\begin{equation*}\label{tp5}\begin{aligned}
&f_2(s,B)/A^4_2=Bc_2=B(B-s^2)^2(2\psi\psi_{ss}-\psi_{s}^2)- B(B-s^2)(6s\psi\psi_{s}+\psi_{ss})+2Bs\psi_{s},\\
&f_2(s,B)=-3B\{9s^4-6(1+2B)s^3+6(1+2B)s^2+2(2B^2-10B-1)s+3B(B+2)\},
\end{aligned}\end{equation*}
$(\cdots)$, independent of $r_{00}^2$, is a polynomial of $A_1, A_2$(or $s, B$), and the degree of $(\cdots)$ in $s$ is no more than deg$(A_1^{j_2}A_2^{k_2})=j_2+k_2$.\\

(3) The calculations for $(\overline{R}^i_{\,\,\,j}b^j)_{.\,i}$.

By \eqref{tc1}, $(\overline{R}^i_{\,\,\,j}b^j)_{.\,i}$ can be rewritten as
\begin{equation}\label{tp6}\begin{aligned}
(\overline{R}^i_{\,\,\,j}b^j)_{.\,i}&=({^\alpha R}^i_{\,\,\,j}b^j)_{.\,i}+\frac{f_3(s,B)}{\alpha^{3}A^5_2}r^2_{00}+\sum_{i_3,j_3,k_3}\frac 1{\alpha^{i_3}A_1^{j_3} A_2^{k_3} }(\cdots),
\end{aligned}\end{equation}
where the sum indices $i_3,j_3$ and $k_3$ satisfy $-1\leq i_3\leq 2, 0\leq j_3\leq 4$ and $0\leq k_3\leq 5$ respectively, $f_3(s,B)$ is a nonzero polynomial of $s, B$ given by
\begin{equation*}\label{tp7}\begin{aligned}
f_3(s,B)&/A^5_2=2(B-s^2)^3\psi\psi_{sss}-(B-s^2)^2(18s\psi\psi_{ss}+6\psi\psi_{s}+\psi_{sss})\\
&+(B-s^2)(5s\psi_{ss}+\psi_{s}+24s^2\psi\psi_{s})-2s^2\psi_{s},\\
f_3(s,B)&=-162s^6+27(16B+5)s^5-45(8B^2+2B-1)s^4-9(4B^2+31B+1)s^3\\
&+9(40B^3-12B^2-9B-1)s^2+9B(-20B^2+58B+7)s-3B(16B^3+12B^2+54B-1),\\
\end{aligned}\end{equation*}
$(\cdots)$, independent of $r_{00}^2$, is a polynomial of $A_1, A_2$(or $s, B$), and the degree of $(\cdots)$ in $s$ is no more than deg$(A_1^{j_3}A_2^{k_3})+1={j_3}+{k_3}+1$.

(4) The calculations for $(\overline{R}^m_{\,\,\,m}b^i)_{.\,i}$.

By \eqref{td1} and a direct computation, we obtain
\begin{equation}\label{tp8}\begin{aligned}
(\overline{R}^m_{\,\,\,m})_{.\,i}b^i=&({^\alpha R}^m_{\,\,\,m})_{.\,i}b^i
+\frac{r_{00}^2}{\alpha^3}\{(B-s^2)c_{2s}-2sc_2\}+\frac{1}{\alpha^2}\{r_{00}s_0((B-s^2)c_{4s}-sc_4)\\
&+r_{00}r_0((B-s^2)c_{6s}-sc_6+2c_2)+r_{00|0}((B-s^2)c_{8s}-sc_8)\}\\
&+\frac{1}{\alpha}\{s^2_0(B-s^2)c_{10s}+rr_{00}((B-s^2)c_{11s}+c_6)+r^2_{0}(-(B-s^2)c_{11s}+2c_6)\\
&+r_0s_{0}((B-s^2)c_{13s}+2c_4)+r_{00}r^m_{\,\,\,m}(B-s^2)c_{14s}+r_{0m}s^m_{\,\,\,0}(-(B-s^2)c_{14s}+c_{16s})\\
&+r_{00|m}b^m((B-s^2)c_{14s}+c_8)+r_{0m|0}b^m(-(B-s^2)c_{14s}+2c_8)\\
&+s_{0|0}c_{18s}+s_{0m}s^m_{\,\,\,0}c_{19s}\}\\
&+rs_0((B-s^2)c_{20s}+sc_{20}+c_{13})+s_{m}s^m_{\,\,\,0}((B-s^2)c_{22s}+sc_{22}+2c_{19}) \\
&+s_{m}r^m_{\,\,\,0}(3(B-s^2)c_{23s}+3sc_{23}-c_{16}) -s_{0}r^m_{\,\,\,m}(2(B-s^2)c_{23s}+2sc_{23}) \\
&+r_{m}s^m_{\,\,\,0}(2(B-s^2)c_{23s}+2sc_{23}+c_{16})-s_{0|m}b^m(2(B-s^2)c_{23s}+2sc_{23}-c_{18}) \\
&+s_{m|0}b^m((B-s^2)c_{23s}+sc_{23}+c_{18})+s^m_{\,\,\,0|m}((B-s^2)c_{24s}+sc_{24})\\
&+(2r_0r^m_{\,\,\,m}-2r_mr^m_{\,\,\,0}+b^kb^m(r_{0k|m}-r_{mk|0}))c_{14}\\
&+\alpha\{s_ms^m((B-s^2)c_{25s}+2sc_{25}-c_{22})+s^i_{\,\,\,m}s^m_{\,\,\,i}((B-s^2)c_{26s}+2sc_{26})+b^ks^m_{\,\,\,k|m}c_{24}\}.  \end{aligned}\end{equation}
Thus, we have
\begin{equation}\label{tp9}\begin{aligned}(\overline{R}^m_{\,\,\,m}b^i)_{.\,i}&=&
({^\alpha R}^m_{\,\,\,m}b^i)_{.\,i}+\frac{f_4(s,B)}{\alpha^3A^5_2}r_{00}^2+\sum_{{i_4},{j_4},{k_4}}\frac 1{\alpha^{i_4}A_1^{j_4}A_2^{k_4}}(\cdots),\end{aligned}\end{equation}
where the sum indices ${i_4},j_4$ and ${k_4}$ satisfy $-1\leq i_4\leq 2, 0\leq j_4\leq 4$ and $0\leq k_4\leq 5$ respectively, $f_4(s,B)$ is a nonzero polynomial of $s, B$ given by
\begin{equation*}\label{tp10}\begin{aligned}
f_4(s,B)&/A^5_2=(B-s^2)c_{2s}-2sc_2\\
&=2(B-s^2)^3\psi\psi_{sss}-(B-s^2)^2(18s\psi\psi_{ss}+6\psi\psi_{s}+\psi_{sss})\\
&+(B-s^2)(6s\psi_{ss}+2\psi_{s}+24s^2\psi\psi_{s})-4s^2\psi_{s},\\
f_4(s,B)&=-162s^6+216(1+2B)s^5-90(4B^2+4B+1)s^4+54(4B^2+B+1)s^3\\
&+18(16B^3-15B^2-9B-1)s^2+54B(-4B^2+9B+1)s-6B(4B^3+24B-1), \nonumber\\
 \end{aligned}\end{equation*}
$(\cdots)$, independent of $r_{00}^2$, is a polynomial of $A_1, A_2$(or $s, B$), and the degree of $(\cdots)$ in $s$ is no more than deg$(A_1^{j_4}A_2^{k_4})+1={j_4}+{k_4}+1$.

Plugging \eqref{tp2}, \eqref{tp4}, \eqref{tp6} and \eqref{tp9} into \eqref{te7} yields
\begin{equation}\label{tp11}
\begin{aligned}
0&={^\alpha R}^i_{\,\,\,j}b^jb_i+\frac{f_1(s,B)}{\alpha^{2}A^4_2}r^2_{00}+\sum_{i_1,j_1,k_1}\frac 1{\alpha^{i_1}A_1^{j_1}A_2^{k_1}}(\cdots)\\
&-\frac{1}{n-1}\{B{^\alpha R}^m_{\,\,\,m}+\frac{f_2(s,B)}{\alpha^{2}A^4_2}r^2_{00}+\sum_{i_2,j_2,k_2}\frac 1{\alpha^{i_2}A_1^{j_2}A_2^{k_2}}(\cdots)\}\\
&-\frac{1}{n+1}\{s\alpha({^\alpha R}^i_{\,\,\,j}b^j)_{.\,i}+\frac{sf_3(s,B)}{\alpha^{2}A^5_2}r^2_{00}+s\sum_{i_3,j_3,k_3}\frac 1{\alpha^{i_3-1}A_1^{j_3} A_2^{k_3} }(\cdots)\}\\
&+\frac{1}{n^2-1}\{s\alpha({^\alpha R}^m_{\,\,\,m}b^i)_{.\,i}+\frac{sf_4(s,B)}{\alpha^2A^5_2}r_{00}^2+s\sum_{{i_4},{j_4},{k_4}}\frac 1{\alpha^{i_4-1}A_1^{j_4}A_2^{k_4}}(\cdots)\}.\\
 \end{aligned}\end{equation}

Multiplying both sides of \eqref{tp11} with $\alpha^2$ yields
\begin{equation}\label{tp12}
\begin{aligned}
0&=\alpha^2\{{^\alpha R}^i_{\,\,\,j}b^jb_i-\frac{1}{n-1}B{^\alpha R}^m_{\,\,\,m}-\frac{1}{n+1}s\alpha({^\alpha R}^i_{\,\,\,j}b^j)_{.\,i}+\frac{1}{n^2-1}s\alpha({^\alpha R}^m_{\,\,\,m}b^i)_{.\,i}\}\\
&+\alpha^2\sum_{i,j,k}\frac 1{\alpha^{i}A_1^{j}A_2^{k}}(\cdots)+r^2_{00}\frac{1}{A^5_2}\{A_2f_1(s,B) -\frac{A_2f_2(s,B)}{n-1}
-\frac{sf_3(s,B)}{n+1}  +\frac{sf_4(s,B)}{n^2-1} \}.
 \end{aligned}\end{equation}
Consequently, the term
\begin{eqnarray}r^2_{00}\frac{1}{A^5_2}\{A_2f_1(s,B) -\frac{A_2f_2(s,B)}{n-1}
-\frac{sf_3(s,B)}{n+1}  +\frac{sf_4(s,B)}{n^2-1} \} \label{tp13}\end{eqnarray}
must be divided by $\alpha^2$. Observe that the coefficient of $r_{00}^2$ in \eqref{tp13} can not be divided by $\alpha$ from the definition of $A_2$ and $f_i(s,b)$ (i=1,...,4). Thus $r_{00}^2$ must be divided by $\alpha$, which means $r_{00}=\sigma\alpha^2$ for some function $\sigma=\sigma(x)$.
 This completes the proof of Lemma \ref{conf}.
 \hfill$\Box$

  \bigskip
\noindent{\it Proof of Theorem \ref{tta1}.} \\
Assume that $F$ is of scalar flag curvature. Note that $B<\frac 14$ for a Matsumoto metric $F$.

By \eqref{tp1}, we conclude that there does not exist $A^{-5}_2$ in the term $\overline{R}^i_{\,\,\,j}b^jb_i$. Also, there does not exist $A^{-5}_2$ in the term ${\overline{R}}^m_{\,\,\,m}$ by \eqref{td1}. So we focus on the term $(\overline{R}^i_{\,\,\,j}b^j)_{.\,i}$ and $(\overline{R}^m_{\,\,\,m}b^i)_{.\,i}$.

By \eqref{tc1}, $(\overline{R}^i_{\,\,\,j}b^j)_{.\,i}$ can be rewritten as
\begin{equation}\label{th1}\begin{aligned}
(\overline{R}^i_{\,\,\,j}b^j)_{.\,i}&=2\psi\psi_{sss}(B-s^2)^3 \frac{(r_{00}-2 \alpha Q s_0)^2}{\alpha^3} +[\cdots]\\
&=\frac{324(B-s^2)^3(r_{00}-2 \alpha Q s_0)^2}{\alpha^3A^5_2}+[\cdots],\\
\end{aligned}\end{equation}
where $[\cdots]$ does not contain the factor $A^{-5}_2$.

For the same reason, by \eqref{tp8}, $(\overline{R}^m_{\,\,\,m}b^i)_{.\,i}$ can be rewritten as
\begin{equation}\label{th2}\begin{aligned}
(\overline{R}^m_{\,\,\,m}b^i)_{.\,i}&=2\psi\psi_{sss}(B-s^2)^3 \frac{(r_{00}-2 \alpha Q s_0)^2}{\alpha^3} +[\cdots]\\
&=\frac{324(B-s^2)^3(r_{00}-2 \alpha Q s_0)^2}{\alpha^3A^5_2}+[\cdots],\\
\end{aligned}\end{equation}
where $[\cdots]$ does not contain the factor $A^{-5}_2$.

Thus, \eqref{te7} can also be rewritten as
\begin{equation}\label{th3}\begin{aligned}
0&=-\frac{1}{n+1}(\frac{324(B-s^2)^3(r_{00}-2 \alpha Q s_0)^2}{\alpha^3A^5_2}
- \frac{324(B-s^2)^3(r_{00}-2 \alpha Q s_0)^2}{(n-1)\alpha^3A^5_2} )s \alpha+[\cdots]\\
&=-\frac{324(n-2)s(B-s^2)^3(r_{00}-2 \alpha Q s_0)^2}{(n^2-1)\alpha^2A^5_2}+[\cdots],\\
&=-\frac{324(n-2)s(1+2B+3s)^5(1+2s)(B-s^2)^3(A_1r_{00}-2 \alpha s_0)^2}{(n^2-1)\alpha^2(1-4s^2)((1+2B)^2-9s^2)^5}+[\cdots],
 \end{aligned}\end{equation}
where $[\cdots]$ does not contain the factor $A^{-5}_2$.

By Lemma \ref{conf}, we have
\begin{equation}\label{a-1} A_1r_{00}-2 \alpha s_0=\alpha(\sigma \alpha-2\sigma\beta-2s_0).   \end{equation}
Note
\begin{equation}\label{a0}  s(1+2B+3s)^5(1+2s)= \frac{16(1+2B)(5+4B)}{9}(1+2B+3s)\,\,mod\,\,\{(1+2B)^2\alpha^2-9\beta^2\}.  \end{equation}
For the same reason in discussing \eqref{tq7}, we get that $[(1+2B)\alpha+3\beta](\sigma \alpha-2\sigma\beta-2s_0)^2$ must be divided by $(1+2B)^2\alpha^2-9\beta^2$ from \eqref{th3}, \eqref{a-1} and \eqref{a0}, i.e.,
\begin{equation}\label{a01}
[(1+2B)\alpha+3\beta](\sigma \alpha-2\sigma\beta-2s_0)^2=0\,\,mod\,\,\{(1+2B)^2\alpha^2-9\beta^2\}.\end{equation}\\

\textbf{Case one}: $\sigma=0$. By \eqref{a01}, we have $s_0^2=\rho\{(1+2B)^2\alpha^2-9\beta^2\}$ for some function $\rho=\rho(x)$. Differentiating both sides yields $(1+2B)^2\rho a_{ij}=9\rho b_ib_j+s_is_j$. Since $n\geq3$, we get $\rho=0$ and $s_0=0$. Hence, $\beta$ is a constant Killing form.\\

\textbf{Case two}: $\sigma\neq0$. Let $q:=[(1+2B)\alpha+3\beta](\sigma \alpha-2\sigma\beta-2s_0)^2$. $q$ can be rewritten as $q=\alpha \,\,q_{even}+q_{odd}$, where
\begin{equation}\label{a1}\begin{aligned}
q_{even}:= (1+2B)\sigma ^2 \alpha^2- 8(1-B)\sigma ^2 \beta^2-4\sigma(1-4B) \beta s_0+4(1+2B)s^2_0
\end{aligned}\end{equation}
and
\begin{equation}\label{a2}\begin{aligned}
q_{odd}:= -(1+8B)\sigma ^2 \alpha^2\beta- 4(1+2B)\sigma \alpha^2 s_0+12\sigma ^2 \beta^3+24\sigma\beta^2 s_0+12\beta s_0^2.
\end{aligned}\end{equation}
Hence, $q_{even}$ can be divided by $(1+2B)^2\alpha^2-9\beta^2$.

Meanwhile, we have
\begin{equation}\label{a3}\begin{aligned}
(\sigma \alpha-2\sigma\beta-2s_0)(\sigma \alpha+2\sigma\beta+2s_0)= \sigma ^2\alpha^2-4\sigma ^2 \beta^2-8\sigma\beta s_0 -4s^2_0
\end{aligned}\end{equation}
can be divided by $(1+2B)^2\alpha^2-9\beta^2$ from \eqref{th3}.

$\frac{1}{2\sigma}\{\eqref{a1}+(1+2B)\times\eqref{a3}\}$ yields
\begin{equation*} \sigma(1+2B) \alpha^2-6 \sigma \beta^2  -6 \beta s_0 =0\,\,mod\,\,\{(1+2B)^2\alpha^2-9\beta^2\}.     \end{equation*}
That is
\begin{equation} \label{a4} \sigma(1+2B) \alpha^2-6 \sigma \beta^2  -6 \beta s_0=\varrho \{(1+2B)^2\alpha^2-9\beta^2\}   \end{equation}
holds for some function $\varrho=\varrho(x)$. \eqref{a4} can be rewritten as
\begin{equation*} \label{a5} (1+2B)[\sigma -(1+2B)\varrho]\alpha^2=3\beta(-3\varrho \beta+2 \sigma \beta -2 s_0),   \end{equation*}
which is equivalent to
\begin{equation} \label{a6} \sigma=(1+2B)\varrho   \end{equation}
and
\begin{equation} \label{a7} 3\varrho \beta=2 \sigma \beta -2 s_0.  \end{equation}
Differentiating both sides of \eqref{a7} and contracting it with $b^i$ yield $3\varrho=2 \sigma$, which is a contraction with \eqref{a6}.

Above all, we have $r_{00}=s_0=0$, i.e., $\beta$ is a constant Killing form. By Lemma \ref{ttq1}, $F$ is projectively flat.

The converse is obvious. This completes the proof of Theorem \ref{tta1}. \hfill$\Box$

\vspace{4mm}
From Theorem \ref{tta1}, one obtains the following
\begin{cor}\label{xxe2}
(\cite{ZS2}) Let $F=\frac{\alpha^2}{\alpha-\beta}$ be a non-Riemannian Matsumoto metric on an $n$-dimensional manifold $M$, $n\geq 3$. Then the following conditions are equivalent:\\
(a) $F$ is of weakly isotropic flag curvature $K$; (b) $F$ is of constant flag curvature $K$; (c) $\alpha$ is a flat metric and $\beta$ is a constant one form.\\
In this case, $F$ is locally Minkowskian.
\end{cor}

{\bf{Acknowledgements}} Author would like to express her sincere thanks to Prof. Yibing Shen for his constant encouragement and many helpful discussions on this paper.

\noindent\author{{ \small Xiaoling Zhang }\\
{\small {\it College of Mathematics and Systems Science, Xinjiang University}}\\
 {\small {\it Urumqi, Xinjiang Province, 830046, P.R.China}}\\
{\small{\it E-mail: xlzhang@ymail.com}}}

\end{document}